\documentclass[letterpaper, 10 pt, conference]{ieeeconf}  

\IEEEoverridecommandlockouts                              

\usepackage{amsmath} 
\usepackage{amssymb}  

\usepackage{amsthm}
\usepackage{nccmath}
\usepackage{mathtools}

\usepackage[noadjust]{cite} 

\usepackage{arydshln} 
\usepackage{mdframed} 
\usepackage{float} 
\newfloat{algorithm}{t}{}
\usepackage{pgfplots}
\newlength\breite

\makeatletter
\def\endthebibliography{%
	\def\@noitemerr{\@latex@warning{Empty `thebibliography' environment}}%
	\endlist
}
\makeatother

\newcommand{\norm}[1]{\left\lVert#1\right\rVert}
\newcommand\numberthis{\addtocounter{equation}{1}\tag{\theequation}}
\renewcommand{\refeq}[1]{\overset{#1}{=}}
\newcommand{\refleq}[1]{\overset{#1}{\leq}}

\title{\LARGE \bf
Data-driven online convex optimization for control of dynamical systems
}

\author{Marko Nonhoff$^{1}$ and Matthias A. M\"uller$^{1}$
\thanks{$^{1}$Both authors are with the Leibniz University Hannover, Institute of Automatic Control, 30167 Hannover, Germany.
        {\tt\small \{nonhoff,mueller\}@irt.uni-hannover.de}}%
}

\allowdisplaybreaks

\newcommand\copyrighttext{%
	\footnotesize \copyright 2021 IEEE. Personal use of this material is permitted. Permission from IEEE must be obtained for all other uses, in any current or future media, including reprinting/republishing this material for advertising or promotional purposes, creating new collective works, for resale or redistribution to servers or lists, or reuse of any copyrighted component of this work in other works.}
\newcommand\copyrightnotice{%
	\begin{tikzpicture}[remember picture,overlay]
		\node[anchor=south,yshift=10pt] at (current page.south) {\fbox{\parbox{\dimexpr\textwidth-\fboxsep-\fboxrule\relax}{\copyrighttext}}};
	\end{tikzpicture}%
}

\begin{document}
\maketitle
\thispagestyle{empty}
\pagestyle{empty}
\copyrightnotice
\begin{abstract}
We propose a data-driven online convex optimization algorithm for controlling dynamical systems. In particular, the control scheme makes use of an initially measured input-output trajectory and behavioral systems theory which enable it to handle unknown discrete-time linear time-invariant systems as well as a priori unknown time-varying cost functions. Further, only output feedback instead of full state measurements is required for the proposed approach. Analysis of the closed loop's performance reveals that the algorithm achieves sublinear regret if the variation of the cost functions is sublinear. The effectiveness of the proposed algorithm, even in the case of noisy measurements, is illustrated by a simulation example.
\end{abstract}

\section{INTRODUCTION}

Application of online convex optimization (OCO) to solve control problems has recently gained significant interest. The OCO framework is an online variant of classical numerical optimization. In particular, the cost function is allowed to be time-varying and a priori unknown (see \cite{Hazan2016,Simonetto2020} for an overview). More specifically, an OCO algorithm has to choose an action at every time instant. Only then, after the action is chosen, the environment reveals the cost function which results in a cost depending on the algorithm's action. This framework is of interest for controller design because of its ability to handle time-varying and unknown cost functions and low computational complexity. Therefore, OCO has been applied to the control of linear dynamical systems with (\cite{Nonhoff2020,Li2020}) and without constraints (\cite{Nonhoff2019,Li2019,Agarwal2019}) with applications in, e.g., power grids \cite{Picallo2020}, data center cooling \cite{Lazic2018}, and robotics \cite{Zheng2020}. Most relevant to this work are \cite{Hazan2020,Simchowitz2020}, where the cost functions as well as the system to be controlled are unknown. In these works, the authors leverage techniques from system identification to learn a model of the system online. 

In this work, we take a different approach, which relies on a result from behavioral systems theory. The so-called fundamental lemma shows that the Hankel matrix consisting of a single persistently exciting input-output trajectory spans the whole vector space of possible trajectories of an LTI system \cite{Willems2005}. This result has recently been used to solve various control problems, e.g., model predictive control \cite{Coulson2019,Berberich2020,Berberich2020IFAC}, state-feedback design \cite{Berberich2020ACC,DePersis2020,Xue2021}, and output matching \cite{Markovsky2008}. We apply the fundamental lemma to the OCO framework in order to control dynamical systems subject to time-varying cost functions, where neither the system nor the cost functions are known to the algorithm.

Our contributions in this paper are threefold. First, we consider an unknown system by leveraging results from data-driven control. Compared to alternative approaches in the literature, we thereby remove the need of a (set-based) model description and of an online estimation process. To the author's best knowledge, this is the first work that combines the fundamental lemma and OCO to solve this problem. Second, we extend the results from \cite{Nonhoff2019} to the case of output feedback instead of full state measurements. This result requires considerable adjustments in the algorithm design. Third, we derive a novel data-based characterization of all steady states of a discrete-time LTI system solely in terms of measured input-output-trajectories. Our subsequent analysis reveals that the proposed algorithm still enjoys sublinear regret without access to a model of the system.

This paper is organized as follows. After stating the required assumptions and definitions in Section~\ref{sec:Setting}, we present and explain our algorithm in Section~\ref{sec:Algorithm}. Our main result, a regret analysis of the proposed algorithm, is stated in Section~\ref{sec:RegretAnalysis}. The closed-loop performance is investigated in Section~\ref{sec:Simulations} by means of a simulation example. Section~\ref{sec:Conclusion} concludes the paper. 
All proofs are given in the Appendix.

\textit{Notation:} We denote the set of all integers in the interval $[a,b]$ by $\mathbb I_{[a,b]}$. The Moore-Penrose-Pseudoinverse of a matrix $A$ is given by $A^\dagger$. For a matrix $A$ and a vector $x$, we denote the Euclidean norm by $\norm x$, whereas $\norm A$ is the induced matrix norm. The identity matrix of size $n \times n$ is given by $I_n \in \mathbb{R}^{n\times n}$, $1_n\in\mathbb{R}^n$ is the vector of all ones, and $0_n \in \mathbb{R}^n$ is the vector of all zeros. A sequence $z = \{ z_k \}_{k=1}^N$, $z_k \in \mathbb{R}^m$, induces the Hankel matrix 
\setlength{\arraycolsep}{2pt}
\[
H_L(z) = \begin{bmatrix} z_1 & z_2 & \dots & z_{N-L+1} \\ z_2 & z_3 & \dots & z_{N-L+2} \\ \dots &  \dots & \dots & \dots \\ z_L & z_{L+1} & \dots & z_N \end{bmatrix} = 
\begin{bmatrix} H_L^1(z) \\ H_L^2(z) \\ \dots \\ H_L^L(z) \end{bmatrix}.
\]
We denote a matrix containing a subset of rows of $H_L(z)$ by $H_L^{a:b}(z) = \begin{bmatrix} (H_L^a(z))^\top & \dots & (H_L^b(z))^\top \end{bmatrix}^\top$.
We also define $z_{[a,b]} = \begin{bmatrix} z_a^\top & \dots & z_b^\top \end{bmatrix}^\top$. With a slight abuse of notation, we write $z$ for the sequence itself as well as for the stacked vector of all components $z_{[1,N]}$. The shift operator $\sigma$ is defined as $\sigma z := \begin{bmatrix} z_2^\top & \dots & z_N^\top \end{bmatrix}^\top$.
For matrices $A$ and $B$, we denote the Kronecker product by $A \otimes B$. 

\section{SETTING} \label{sec:Setting}

We consider discrete-time LTI systems of the form
\[
	x_{t+1} = Ax_t + Bu_t, \qquad	y_{t} = Cx_t + Du_t, \numberthis \label{eq:LTIsystem}
\]
where $x\in\mathbb{R}^n$, $u\in\mathbb{R}^m$, and $y\in\mathbb{R}^p$. Throughout this paper, we require $(A,B)$ to be controllable and $(A,C)$ to be observable. However, we assume that the matrices $(A,B,C,D)$ in \eqref{eq:LTIsystem} as well as the initial condition $x_0 = \bar x$ are unknown and only measurement data in the form of a prerecorded input-output sequence $\{u^d_k,y^d_k\}_{k=0}^{N-1}$ and an upper bound on the system order $n$ are available. We aim to solve the optimal control problem
\[\min_{u\in\mathbb R^{m(T+1)}} \sum_{t=0}^T L_t(u_t,y_t) \quad \text{s.t. } \eqref{eq:LTIsystem}. \numberthis \label{eq:OCP} \]
The main difficulty arises from the fact that the cost functions $L_t: \mathbb{R}^m \times \mathbb{R}^p \rightarrow \mathbb{R}$ are not known a priori. In specific, at every time $t\in\mathbb{I}_{[0,T]}$, our algorithm has to decide on a control action $u_t$ which is then applied to system \eqref{eq:LTIsystem}. Only then, after $u_t$ is chosen by the algorithm, the output $y_t$ is measured and the cost function $L_t(u,y)$ is revealed resulting in the cost $L_t(u_t,y_t)$.
We denote the optimal input and output sequence in hindsight, i.e., the solution to \eqref{eq:OCP} with full knowledge of the cost functions, by $u^* = \{u^*_t\}_{t=0}^T$ and $y^* = \{y^*_t\}_{t=0}^T$. In order to solve the problem at hand, we require the data sequence to be persistently exciting as captured in the following definition and Assumption~4.

\noindent {\bf Definition 1:} A signal $\{u_k\}_{k=0}^{N-1}$, $u_k\in\mathbb R^m$, is persistently exciting of order $L$ if $\text{rank}(H_L(u)) = mL$.

This definition allows for a characterization of all system trajectories of \eqref{eq:LTIsystem} as given in the following result. It shows that the Hankel matrix induced by a single persistently exciting input sequence spans the vector space of all system trajectories of \eqref{eq:LTIsystem}. The result was first published in the context of behavioral systems theory \cite{Willems2005} and can be formulated in the classical state-space setting as follows.

\noindent {\bf Theorem 1:} (\cite{Berberich2020ECC}, Theorem 3) Suppose $\{u_k^d,y_k^d\}_{k=0}^{N-1}$ is a trajectory of system \eqref{eq:LTIsystem}, where $u^d$ is persistently exciting of order $L+n$. Then, $\{u_k,y_k\}_{k=0}^{L-1}$ is a trajectory of \eqref{eq:LTIsystem} if and only if there exists $\alpha \in \mathbb{R}^{N-L+1}$ such that
\[
	\begin{bmatrix} H_L(u^d) \\ H_L(y^d) \end{bmatrix} \alpha = \begin{bmatrix} u \\ y \end{bmatrix}.
\]

As discussed below, we aim to track a sequence of desired steady states without access to a model of the system. We use the same data-driven definition of steady states as \cite{Berberich2020}.

\noindent {\bf Definition 2:} An input-output pair $(u^s,y^s)$ is an equilibrium of \eqref{eq:LTIsystem}, if the sequence $\{u_k,y_k\}_{k=0}^n$ with $(u_k,y_k) = (u^s,y^s)$ for all $k \in \mathbb{I}_{[0,n]}$ is a trajectory of \eqref{eq:LTIsystem}.

Definition~2 states that an input-output pair $(u^s,y^s)$ is an equilibrium if and only if a sequence consisting of $(u^s,y^s)$ for at least $n+1$ consecutive time steps is a valid trajectory of \eqref{eq:LTIsystem}. Moreover, we require that for every output $y$, there exists at least one corresponding input and internal state such that $y$ is an equilibrium of \eqref{eq:LTIsystem}. This assumption is used in tracking problems, compare, e.g., \cite[Section 1.5]{Rawlings2017}.

\noindent {\bf Assumption 1:} For every $y \in \mathbb{R}^p$, there exists $u\in\mathbb{R}^m$ such that the input-output pair $(u,y)$ is an equilibrium of \eqref{eq:LTIsystem}. 

\textit{Remark 1:} A sufficient condition for this assumption to be met is that the matrix $\begin{bmatrix} I_n-A & -B \\ C & D \end{bmatrix}$ has full row rank. This condition requires $m \geq p$, i.e., there have to be at least as many inputs to the system as outputs \cite[Section 1.5]{Rawlings2017}.

Concerning the cost functions $L_t$, we require the following technical assumptions, which are fairly standard in OCO.

\noindent {\bf Assumption 2:} The functions $L_t$ satisfy $L_t(u,y) = f_t^u(u) + f_t^y(y)$, where $f_t^\xi(\xi)$ is $\alpha_\xi$-strongly convex, $l_\xi$-smooth, and Lipschitz continuous with Lipschitz constant $L_\xi$, $\xi = u,y$.

\textit{Remark 2:} We assume Lipschitz continuity for clarity of exposition of our results although $l$-smoothness and Lipschitz continuity cannot be satisfied \textit{globally} simultaneously. However, under the reasonable assumption that $y_t$ and $u_t$ remain within a bounded set for all times, it suffices if Assumption~2 is satisfied on this set. Moreover, the same techniques as in \cite{Nonhoff2019} can be used to avoid assuming Lipschitz continuity. In this case, all triangle inequalities in the proofs are replaced by Jensen's inequality, which entails additional assumptions on the step sizes and the regularity of the cost functions $l/\alpha$.

Note that $l$-smoothness implies Lipschitz continuity of the gradients $\nabla f_t$ with constant $l$. Let $\theta_t = \arg\min_y f_t^y(y)$ and $\eta_t = \arg\min_u f_t^u(u)$ be the point-wise in time minimizers of the cost function $L_t(u,y)$. We require the input-output pair $(\eta_t,\theta_t)$ to be an equilibrium of \eqref{eq:LTIsystem}. Thus, minimizing \eqref{eq:OCP} amounts to tracking a sequence of a priori unknown steady states of~\eqref{eq:LTIsystem}. Generalizing this assumption to allowing arbitrary cost functions (termed \textit{economic} cost functions in the context of MPC \cite{Rawlings2017}) is part of our ongoing work.

\noindent {\bf Assumption 3:} The input-output pair $(\eta_t,\theta_t)$ is an equilibrium of system \eqref{eq:LTIsystem} for all $t \in \mathbb{I}_{[0,T]}$.

\section{ALGORITHM} \label{sec:Algorithm}

\begin{algorithm} 
	\vspace{4pt}
	{
	\setlength\belowdisplayskip{0pt}
	\setlength\abovedisplayskip{0pt}
	\fbox{\parbox{.95\linewidth}{
	\underline{Algorithm 1: Data-Driven Output Feedback}
	
	\noindent Given stepsizes $\gamma_u,\gamma_y$, data $(u^d,y^d)$ and an initialization $v_0\in\mathbb R^m$, $\hat u_{-1} \in \mathbb{R}^{m(\mu+1)}$,
	\begin{align}
		&v_t = v_{t-1} - \gamma_u \nabla f_{t-1}^u(v_{t-1}) \label{eq:InputOGD} \\
		&\text{Choose $\alpha_t$ and $\omega_t$ such that} \nonumber \\[-4pt]
		&\quad H_\alpha \omega_t = \begin{bmatrix} 0_{mn}^\top & (1_{n+\mu+1} \otimes v_t)^\top & 0_{pn}^\top \end{bmatrix}^\top \label{eq:defomega} \\
		&\quad H_\alpha \alpha_t  = \begin{bmatrix} u_{[t-n:t-1]} \\ \hdashline \sigma \hat u_{t-1} \\ 1_{n+1} \otimes (u^s_{t-1} - v_{t-1}) \\ \hdashline y_{[t-n:t-1]} \end{bmatrix} \label{eq:defalpha} \\
		&\hat y_{t}^\mu = Y^{n+\mu+1} (\alpha_t + \omega_t) \label{eq:OutputPrediction} \\
		&y^s_t = \hat y_{t}^\mu - \gamma_y \nabla f_{t-1}^y(\hat y_{t}^\mu) \label{eq:OutputOGD}\\
		&u^s_t = (I_m - S_u^\dagger S_u) v_t - S_u^\dagger S_y y^s_t \label{eq:SSInput} \\
		&\beta_t \hspace{-3pt} = \hspace{-3pt}
		\begin{cases}
			&\hspace{-10pt} \arg\min_\beta \hspace{5ex}{\norm{Q\beta}} \\
			&\hspace{-10pt}\text{s.t. } H_\beta \beta = \hspace{-2pt} \begin{bmatrix} 0_{mn} \\  1_{n+1} {\otimes} u_t^s {-} U^{\tilde \mu:\tilde \mu+n} (\alpha_t+\omega_t)  \\ 0_{pn} \\  1_{n} {\otimes} y_t^s - Y^{\tilde \mu:\tilde \mu+n-1}(\alpha_t + \omega_t)  \end{bmatrix}
		\end{cases} \label{eq:defbeta} 
		\\
		&\hat u_t = \begin{bmatrix} \sigma \hat u_{t-1} \\ u^s_{t-1} - v_{t-1} \end{bmatrix} + U^{n+1:n+\mu+1}\beta_t \label{eq:PredInputs} \\
		&u_t = \hat u_{t,1} + v_t \label{eq:OutputAlgo}
	\end{align}
	}}
	}
	\vspace{-20pt}
\end{algorithm}

First, let $U = H_{2n+\mu+1}(u^d)$ and $Y = H_{2n+\mu+1}(y^d)$ be the Hankel matrices associated with the input and output data, respectively, $H_\alpha = \begin{bmatrix} U^{1:n} \\ U^{n+1:\tilde \mu +n} \\ Y^{1:n} \end{bmatrix}$, $H_\beta = \begin{bmatrix} U^{1:n} \\ U^{\tilde\mu:\tilde\mu+n} \\ Y^{1:n} \\ Y^{\tilde\mu:\tilde\mu+n-1} \end{bmatrix}$, and $\tilde \mu = n+\mu+1$. The proposed data-driven OCO scheme is given in Algorithm~1. In the framework described above, it computes a control input $u_t$ based on previous cost functions $L_{t-1}(u,y) = f_{t-1}^u(u)+f_{t-1}^y(y)$ and measurements $y_{[t-n,t-1]}$. After $u_t$ is applied to system \eqref{eq:LTIsystem}, a new cost function $L_t$ is revealed, which results in the cost $L_t(u_t,y_t)$. 
Roughly speaking, Algorithm~1 employs online gradient descent (OGD) \cite{Hazan2016} twice to track the optimal equilibrium. In \eqref{eq:InputOGD}, OGD is used to compute an estimate $v_t$ of the optimal input $\eta_t$. Then, $\hat y_t^\mu \in\mathbb{R}^p$ is computed in \eqref{eq:defomega} - \eqref{eq:OutputPrediction}, where $\hat y_t^\mu$ is a prediction $\mu$ time steps ahead and predicted at time $t$. Therefore, $\mu$ is the prediction horizon of Algorithm~1. 
Note that, at every time step $t$, the algorithm computes an initial condition ($n$ time steps), predicts the future output $\hat y_t^\mu$ ($\mu$ time steps) and enforces a terminal constraint ($n+1$ time steps). In order to apply Theorem~1, we thus require persistency of excitation of order $3n+\mu+1$ and a sufficiently long prediction horizon.

\noindent {\bf Assumption 4:} The input $u^d$ of the data sequence is persistently exciting of order $3n+\mu+1$.


\noindent {\bf Assumption 5:} The prediction horizon satisfies $\mu \geq \mu^*$, where $\mu^*$ is the controllability index of the system.

\textit{Remark 3:} Since $\mu^* \leq n$ and we require an upper bound for the order of the system to be known, an upper bound for $\mu^*$ is always available. However, a smaller prediction horizon can be beneficial for the algorithm's performance.

Next, a desired equilibrium is computed in \eqref{eq:OutputOGD} - \eqref{eq:SSInput}. In \eqref{eq:OutputOGD}, the predicted output $\hat y_t^\mu$ is improved using OGD. The resulting output $y^s_t$ is always a steady-state output by Assumption~1. In \eqref{eq:SSInput}, Algorithm~1 determines the corresponding steady-state input as shown by the following result.

\noindent {\bf Lemma 1:} Assume that the sequence $u^d$ is persistently exciting of order $2n+1$. Then, the input-output pair $(u,y)$ is an equilibrium of \eqref{eq:LTIsystem} if and only if
\[
	S_u u + S_y y = 0,
\]
where $ \begin{bmatrix} S_u & S_y \end{bmatrix} = \left( H_{n+1}H_{n+1}^\dagger - I_{(m+p)(n+1)}\right) \begin{bmatrix} \hat I_m & 0 \\ 0 & \hat I_p \end{bmatrix}$, $H_{n+1} = \begin{bmatrix} H_{n+1}(u^d) \\ H_{n+1}(y^d) \end{bmatrix}$, $\hat I_m = 1_{n+1} \otimes I_m$, and $\hat I_p = 1_{n+1} \otimes I_p$. Moreover, let $y\in\mathbb{R}^p$ be any steady-state output of \eqref{eq:LTIsystem} and $v \in \mathbb{R}^m$. The solution to
\[
\min_{u} \norm{u - v} \hspace{4ex} \text{s.t. } (u, y) \text{ is a steady state}
\]
is
\vspace{-1.5ex}
\[
	u = (I_m - S_u^\dagger S_u) v - S_u^\dagger S_y y. \numberthis \label{eq:Lemma1}
\]

The proof is given in the appendix.

\textit{Remark 4:} To the best of the authors' knowledge, the first result in Lemma~1 is a novel data-driven characterization of all steady states of any controllable discrete-time LTI system. This result may be useful in other applications than OCO as well. Moreover, the second result in Lemma~1 implies in particular $\eta_t = (I_m - S_u^\dagger S_u) \eta_t - S_u^\dagger S_y \theta_t$ for $v=\eta_t$ since $(\eta_t,\theta_t)$ is an equilibrium of \eqref{eq:LTIsystem} by Assumption~3.

In \eqref{eq:defbeta}, Algorithm~1 modifies the input sequences used for prediction in \eqref{eq:defomega} - \eqref{eq:defalpha} such that system \eqref{eq:LTIsystem} reaches the desired equilibrium at the end of the prediction horizon, i.e., the sequences given by $\begin{bmatrix} U^\top & Y^\top \end{bmatrix}^\top (\alpha_t+\beta_t+\omega_t)$ reach $(u^s_t,y^s_t)$ in $\mu$ steps. Note that \eqref{eq:defbeta} is always feasible: Equations \eqref{eq:defomega} and \eqref{eq:defalpha} have a solution by Theorem~1, because only an initial condition and the input sequence for the next $\tilde \mu$ time steps are specified by the respective right-hand sides (compare, e.g., the simulation problem in \cite{Markovsky2008}). Moreover, $\begin{bmatrix} U^\top & Y^\top \end{bmatrix}^\top (\alpha_t+\beta_t+\omega_t)$ specifies a sequence initialized by $u_{[t-n,t-1]}$, $y_{[t-n,t-1]}$, which remains at the equilibrium $(u^s_t,y^s_t)$ for the last $n+1$ time steps. This is a valid system trajectory due to controllability and Assumption~5. Since sums of trajectories of LTI systems are valid trajectories themselves, the input-output sequence specified in \eqref{eq:defbeta} is a trajectory of \eqref{eq:LTIsystem}. Therefore, the optimization problem in \eqref{eq:defbeta} always has a solution by Theorem~1. Lastly, in \eqref{eq:PredInputs} an input sequence $\hat u_t \in \mathbb{R}^{m(\mu+1)}$ is collected which is used to predict $\hat y^\mu_{t+1}$ at the next time step. In \eqref{eq:OutputAlgo}, the first part of this input sequence $\hat u_{t,1}$ is applied to system \eqref{eq:LTIsystem}. 

In the following, we derive explicit solution formulas for \eqref{eq:defomega}, \eqref{eq:defalpha}, and \eqref{eq:defbeta}. As discussed earlier, \eqref{eq:defomega}, \eqref{eq:defalpha}, and \eqref{eq:defbeta} always have a solution. One such solution of \eqref{eq:defomega} and \eqref{eq:defalpha} is given by
\[
	\omega_t = H_\alpha^\dagger \begin{bmatrix} 0_{mn} \\ 1_{\tilde\mu} \otimes v_t \\ 0_{pn} \end{bmatrix}, \alpha_t = H_\alpha^\dagger \begin{bmatrix} u_{[t-n:t-1]} \\ \hdashline \sigma \hat u_{t-1} \\ 1_{n+1} \otimes (u^s_{t-1} - v_{t-1}) \\ \hdashline y_{[t-n:t-1]} \end{bmatrix}.
\]
Finally,  \eqref{eq:defbeta} can be solved by the weighted pseudoinverse \cite{Elden1982} by defining $g_t = \begin{bmatrix} 0_{mn} \\ 1_{n+1} \otimes u_t^s - U^{\tilde \mu:\tilde \mu+n} (\alpha_t+\omega_t) \\ 0_{pn} \\ 1_{n} \otimes y_t^s - Y^{\tilde \mu:\tilde \mu+n-1}(\alpha_t + \omega_t) \end{bmatrix}$ by
\[
	\beta_t {=} \hspace{-2pt} \left( I_{N{-}2n{-}\mu} {-} \left( Q \left( I_{N-2n-\mu} {-} H_\beta^\dagger \hspace{-1pt} H_\beta \right) \hspace{-1pt} \right)^\dagger \hspace{-4pt} Q \right) \hspace{-2pt} H_\beta^\dagger g_t. \numberthis \label{eq:PINVbeta}
\]

Since all required pseudoinverses depend only on a priori data, they can be computed offline. Therefore, the calculations that have to be carried out online reduce to two gradient evaluations and multiple matrix-vector multiplications.

\section{REGRET ANALYSIS} \label{sec:RegretAnalysis}

In this section, we analyze the closed-loop performance of Algorithm~1. We define the dynamic regret $\mathcal R$ as a measure of the closed loop's performance by $\mathcal{R} := \sum_{t=0}^T L_t(u_t,y_t) - L_t(u_t^*,y_t^*)$. The dynamic regret is a measure of the performance lost due to not knowing the cost functions $L_t$ a priori. This definition of dynamic regret is in line with the definition in, e.g., \cite{Nonhoff2019,Li2019}. Our main result shows that Algorithm~1 achieves a regret bound which depends linearly on $\sum_{t=0}^T \norm{\eta_t - \eta_{t-1}}$ and $\sum_{t=0}^T \norm{\theta_t - \theta_{t-1}}$. These quantities can be seen as a measure of the variation of the cost functions. Therefore, Algorithm~1 achieves sublinear dynamic regret if this variation is sublinear.

\noindent \textbf{Theorem 2:} Let Assumptions 1-5 be satisfied and choose $\gamma_u \leq \frac{2}{l_u+\alpha_u}$ and $\gamma_y \leq \frac{2}{l_y+\alpha_y}$. Algorithm~1 achieves
\[
\mathcal R \leq C_\mu + C_u \sum_{t=0}^T \norm{\eta_t - \eta_{t-1}} + C_y \sum_{t=0}^T \norm{\theta_t - \theta_{t-1}}
\]
where $C_\mu, C_u,C_y < \infty $ are constants independent of $T$.

The proof is given in the appendix. Theorem~2 is well aligned with other results on dynamic regret in the literature, compare, e.g., \cite{Nonhoff2019,Li2019}. Note that this result implies convergence to the optimal equilibrium if the minima are constant, i.e., $(\eta_{t},\theta_{t}) = (\eta_{t'},\theta_{t'})$ for some $t' \in \mathbb I$ and all $t\geq t'$.

\section{SIMULATIONS} \label{sec:Simulations}

\begin{figure}
	\centering \small
	\setlength\breite{.4\textwidth}
	\vspace{5pt}
	\input{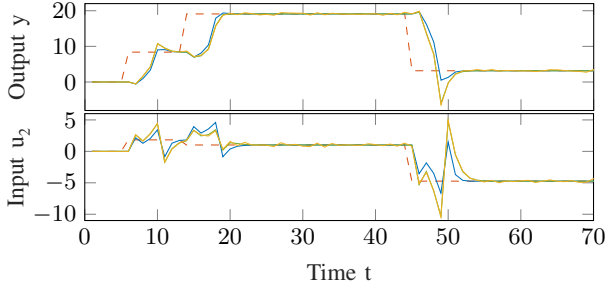}
	\vspace{-55pt}
	\caption{Optimal outputs $\theta_t$ and inputs $\eta_t$ (red dashed) and closed-loop output (top figure) and input trajectory (bottom figure) when applying Algorithm~1 without noise (blue solid), with noisy Hankel matrices (green), and with noisy measurements (yellow).}
	\label{fig:SimResults}
	\vspace{-15pt}
\end{figure}

In this section, we demonstrate the performance of the proposed algorithm and study the effect of noisy measurements by means of simulations. The controlled system is chosen randomly with fixed matrices $(B,C,D)$ with $m=2$, $p=1$ and $A \in \mathbb R^{5\times5}$ is generated by sampling each entry from a uniform distribution over the interval $[-1,1]$. In this case, we get an unstable matrix with $\norm{A} \approx 2.5$. We choose $\gamma_y = \gamma_u = 0.75$, $\mu=5$, $Q=\begin{bmatrix}100 U^{n+1:n+\mu} \\ 100 Y^{n+1:n+\mu} \\ I_{N-2n-\mu+1}\end{bmatrix}$, and $N = 100$. The first two entries in $Q$ penalize the output and the input sequence generated by $\beta_t$ in the transient phase, i.e., before reaching $u^s_t$ and $y_t^s$, respectively. The identity matrix can be seen as a regularization which improves robustness of Algorithm~1 with respect to noise. The cost functions are chosen as $f_t^y(y) = \frac{1}{2} \norm{y-\theta_t}^2$ and $f_t^u(u) = \frac{1}{2} \norm{u-\eta_t}^2$. The optimal input sequence $\eta$ is generated randomly and the optimal output sequence $\theta$ is chosen such that Assumption~3 is satisfied. We consider the following cases: 1)~The algorithm has access to perfect (noiseless) data $y^d$ and online measurements $y_t$, 2)~The output data used to construct the Hankel matrices is noisy, i.e., only $\tilde y^d = y^d + \epsilon^d$ is available, where the random sequence $\epsilon^d = \{\epsilon_k^d\}_{k=0}^{N-1}$ satisfies $\norm{\epsilon^d_k} \leq 10^{-5}$, and 3) Only noisy data $\tilde y^d$ and noisy measurements $\tilde y_t = y_t + \epsilon_t$ are accessible, where $\epsilon_t\in\mathbb{R}$ is generated randomly at each time step and satisfies $\norm{\epsilon_t}\leq 10^{-2}$. The simulation results are shown in Figure~\ref{fig:SimResults}. The algorithm tracks the optimal equilibrium $(\eta_t,\theta_t)$ whenever the cost functions change in all simulations. The closed-loop performance deteriorates when only a noisy data sequence is available, whereas the (additional) effect of measurement noise on the closed loop is small.

\section{CONCLUSIONS} \label{sec:Conclusion}

In this paper, we proposed a novel data-driven OCO scheme for controlling linear dynamical systems. Only a single persistently exciting data trajectory and an upper bound on the system dimension instead of a model of the system and only output instead of full state measurements are required. Our algorithm achieves a similar sublinear regret bound as a model-based algorithm. Future work includes obtaining theoretical guarantees in case of noisy data and relaxing Assumption~3 to allow economic cost functions.



\section*{APPENDIX}

\subsection{Proof of Lemma 1}
By Definition~2 and Theorem~1, $(u,y)$ is a steady state if and only if there exists a $\nu \in \mathbb{R}^{N-n}$ which satisfies
\[
H_{n+1} \nu = \begin{bmatrix} H_{n+1}(u^d) \\ H_{n+1}(y^d) \end{bmatrix} \nu = \begin{bmatrix} \hat I_m u \\ \hat I_p y \end{bmatrix}. \numberthis \label{eq:defSS}
\]
The general solution to this equation is given by
\[
\nu = H_{n+1}^\dagger \begin{bmatrix} \hat I_m u \\ \hat I_p y \end{bmatrix} + \left(I_{N-n} - H_{n+1}^\dagger H_{n+1} \right) \nu'
\]
for an arbitrary $\nu'$. Inserting this solution into \eqref{eq:defSS} yields
\[
	\left(H_{n+1} H_{n+1}^\dagger - I_{(m+p)(n+1)} \right) \begin{bmatrix} \hat I_m & 0 \\ 0 & \hat I_p \end{bmatrix} \begin{bmatrix} u \\ y \end{bmatrix} = 0,
\]
which proves the first result. The optimization problem in the second part of Lemma~1 can therefore be rewritten as $\{\min_{u} \norm{u - v} \text{ s.t. } S_u u + S_y y = 0\}$, which is solved by substituting $\bar z = u -v$ and the fact that $\{\min_x \norm{x} \text{ s.t. } Ax=b\}$ is solved (when being feasible) by $x = A^\dagger b$. \qed

\subsection{Proof of Theorem 2}

We shorten notation by defining $\Theta_\tau := \sum_{t=0}^\tau \norm{\theta_t-\theta_{t-1}}$ and $N_\tau := \sum_{t=0}^\tau \norm{\eta_t-\eta_{t-1}}$. We use the following key result on the convergence rate of gradient descent from \cite{Nesterov2018} to derive a regret bound for Algorithm~1. For an $\alpha$-convex, $l$-smooth function $f: \mathbb R^n \rightarrow \mathbb R$ to be minimized, one gradient descent step $z_1 = z_0 - \gamma \nabla f(z_0)$, where $\gamma \leq \frac{2}{\alpha+l}$, satisfies
\begin{align}
	\norm{z_1 - \theta} \leq \kappa \norm{z_0 - z^*}, \label{eq:ContractionGD}
\end{align}
where $z^* = \arg \min_{z} f(z)$ and $\kappa = 1-\alpha \gamma$. Accordingly, we define $\kappa_y = 1-\alpha_y \gamma_y$ and $\kappa_u = 1-\alpha_u \gamma_u$. By optimality of $(\eta_t,\theta_t)$ and Lipschitz continuity, it holds that
\begin{align*}
	\mathcal R	&\leq C_\mu + L_u \sum_{t=0}^T \norm{u_t-\eta_t} + L_y \sum_{t=0}^{T-\mu} \norm{y_{t+\mu} - \theta_{t+\mu}},
\end{align*}
where $C_\mu = L_y \sum_{t=0}^{\mu-1} \norm{y_t - \theta_t}$ only depends on the initialization and $\{\eta_t,\theta_t\}_{t=0}^{\mu-1}$. Crucially, $C_\mu < \infty$ is a constant independent of $T$. Applying the triangle inequality yields
\begin{align*}
	\mathcal R &\leq C_\mu + L_u \sum_{t=0}^T \norm{u_t-\eta_t} + L_y \sum_{t=0}^{T-\mu} \norm{y_{t+\mu} - \hat y_{t}^\mu} \\
	&+ L_y \sum_{t=0}^{T-\mu} \norm{\hat y_{t}^\mu - \theta_{t-1}} + L_y \sum_{t=0}^{T-\mu} \norm{\theta_{t+\mu}-\theta_{t-1}}. \numberthis \label{eq:RegretUnfinished}
\end{align*}
Using the triangle inequality again we obtain
\[
	\sum_{t=0}^{T-\mu} \norm{\theta_{t+\mu} - \theta_{t-1}} \leq (\mu + 1) \Theta_{T-\mu}. \numberthis \label{eq:thetamu-theta}
\]
In the following, we bound the remaining three sums in the regret bound~\eqref{eq:RegretUnfinished}. First, we analyze the regret of the control inputs $\norm{u_t {-} \eta_t}$. In the second part, we derive a bound on the predicted regret $\norm{\hat y_t^\mu {-} \theta_{t{-}1}}$. Third, we bound the deviation of the predicted output $\hat y^\mu_t$ from the real output $y_{t{+}\mu}$.

\subsubsection*{Part 1 (Input Regret)}
Before we analyze the control inputs, we first derive some useful auxiliary results. We have
\[
	\sum_{t=0}^\tau \norm{v_{t} - \eta_t} \refleq{\eqref{eq:InputOGD},\eqref{eq:ContractionGD}} \kappa_u \sum_{t=0}^\tau \norm{v_{t-1} - \eta_{t-1}} + \sum_{t=0}^\tau \norm{\eta_t - \eta_{t-1}}.
\]
Since $\eta_{-1}$ is not defined yet, we choose $\eta_{-1} = v_{-1}$. Due to $1-\kappa_u>0$, we obtain
\[
	\sum_{t=0}^\tau \norm{v_{t} - \eta_t} \leq \frac{1}{1-\kappa_u} N_\tau. \numberthis \label{eq:vt-etat}
\]
Moreover, by optimality of $\eta_t$ we have $\nabla f_t^u(\eta_t)=0$. Hence, Lipschitz continuity of the gradients yields
\[
	\sum_{t=0}^\tau \norm{v_{t+1} {-} v_t}  \refeq{\eqref{eq:InputOGD}}  \sum_{t=0}^\tau  \norm{\nabla f_t^u(v_t) {-} \nabla f_t^u(\eta_t)} \refleq{\eqref{eq:vt-etat}}  \frac{l_u}{1{-}\kappa_u}  N_\tau. \numberthis \label{eq:vt-vt}
\]
This result immediately implies
\begin{align*}
	&\sum_{t=0}^\tau \sum_{i=0}^\mu \norm{v_{t+i}-v_t} 
	\leq \sum_{t=0}^{\tau} \sum_{i=1}^\mu \sum_{j=1}^i \norm{v_{t+j} - v_{t+j-1}} \\
	&\leq \mu^2 \sum_{t=0}^{\tau+\mu} \norm{v_{t+1} - v_{t}} \refleq{\eqref{eq:vt-vt}} \frac{\mu^2 l_u}{1-\kappa_u} N_{\tau+\mu}. \numberthis \label{eq:vti-vt}
\end{align*}
Next, we show 
\[Y^{\tilde\mu:\tilde\mu+n-1}(\alpha_t + \omega_{t-1}) = Y^{\tilde\mu+1:\tilde\mu+n}(\alpha_{t-1}+\beta_{t-1}+\omega_{t-1}).\]
Since $u_{t-1} \refeq{\eqref{eq:defomega},\eqref{eq:defalpha},\eqref{eq:PredInputs},\eqref{eq:OutputAlgo}} U^{n+1} (\alpha_{t-1} + \beta_{t-1} + \omega_{t-1})$,
we have that $\begin{bmatrix} U^{1:n} \\ Y^{1:n} \end{bmatrix} (\alpha_t + \omega_{t-1}) \refeq{\eqref{eq:defomega},\eqref{eq:defalpha}} \begin{bmatrix} u_{[t-n,t-1]} \\ y_{[t-n,t-1]} \end{bmatrix}$ and
\[
	\begin{bmatrix} U^{2:n+1} \\ Y^{2:n+1} \end{bmatrix} (\alpha_{t-1} + \beta_{t-1} + \omega_{t-1}) \refeq{\eqref{eq:defomega},\eqref{eq:defalpha},\eqref{eq:defbeta}} \begin{bmatrix} u_{[t-n,t-2]} \\ u_{t-1} \\ \hdashline y_{[t-n,t-2]} \\y_{t-1} \end{bmatrix}.
\]
Therefore, the initial conditions of both sequences are the same. Additionally, the input sequences satisfy
\begin{align*}
	&U^{n+1:2n+\mu}(\alpha_t + \omega_{t-1}) \\ 
	\refeq{\eqref{eq:defomega},\eqref{eq:defalpha}} 
	&\begin{bmatrix} \sigma \hat u_{t-1} + U^{n+1:\tilde\mu-1} \omega_{t-1} \\ 1_{n} \otimes (u_{t-1}^s - v_{t-1}) + 1_{n} \otimes v_{t-1} \end{bmatrix} \\
	\refeq{\eqref{eq:defalpha},\eqref{eq:defbeta},\eqref{eq:PredInputs}} &\begin{bmatrix} U^{n+2:\tilde\mu} (\alpha_{t-1} + \beta_{t-1} + \omega_{t-1})  \\ U^{\tilde\mu+1:\tilde\mu+n} \left(\alpha_{t-1} + \beta_{t-1} + \omega_{t-1} \right) \end{bmatrix},
\end{align*}
which shows that the input sequences are the same as well. Therefore, the resulting output sequence must be identical. Thus, $Y^{\tilde\mu:\tilde\mu+n-1}(\alpha_t + \omega_{t-1}) = Y^{\tilde\mu+1:\tilde\mu+n}(\alpha_{t-1}+\beta_{t-1}+\omega_{t-1})$
which implies 
\begin{align}
	&Y^{\tilde\mu:\tilde\mu+n-1}(\alpha_t + \omega_t) - Y^{\tilde\mu:\tilde\mu+n-1}(\omega_t - \omega_{t-1}) \nonumber \\
	= &Y^{\tilde\mu+1:\tilde\mu+n}(\alpha_{t-1}+\beta_{t-1}+\omega_{t-1}). \label{eq:RecursivePredictions}
\end{align}
Moreover, due to the fact that $\begin{bmatrix} U^\top & Y^\top \end{bmatrix}^\top (\alpha_t+\beta_t+\omega_t)$ generates a sequence which is at equilibrium for the last $n+1$ time steps, we obtain
\[
	Y^{\tilde\mu:\tilde\mu+n}(\alpha_t + \beta_t + \omega_t) = 1_{n+1} \otimes y^s_t. \numberthis \label{eq:TerminalOutput}
\]
Next, we bound the difference of the steady-state outputs at consecutive time steps. Defining $y^s_0 = \hat y^\mu_0 = \theta_{-1}$ yields
\begin{align*}
	&\sum_{t=0}^\tau \norm{y^s_{t} -\theta_{t}} 
	\refleq{\eqref{eq:OutputOGD},\eqref{eq:ContractionGD}} \kappa_y \sum_{t=0}^\tau \norm{\hat y^\mu_{t+1} - \theta_{t}} + \sum_{t=0}^\tau \norm{\theta_{t} -\theta_{t-1}}  \\
	&\refleq{\eqref{eq:OutputPrediction},\eqref{eq:RecursivePredictions},\eqref{eq:TerminalOutput}} \kappa_y \sum_{t=0}^\tau \norm{y_{t}^s + Y^{\tilde\mu}(\omega_{t+1} - \omega_{t}) - \theta_{t}} + \Theta_\tau \\
	&\refleq{\eqref{eq:vt-vt}} \quad \kappa_y \sum_{t=0}^\tau \norm{y_t^s - \theta_{t}} + \Theta_\tau + \frac{c_1\kappa_y \sqrt{\tilde\mu}}{2(1-\kappa_u)} N_\tau.
\end{align*}
By rearranging and defining $c_1 := 2 \norm{Y^{\tilde\mu}H_\alpha^\dagger} l_u $ we get
\[
\sum_{t=0}^\tau \norm{y^s_{t} -\theta_{t}} \refleq{\eqref{eq:vt-vt}} \frac{c_1\kappa_y\sqrt{\tilde\mu}}{2(1-\kappa_y)(1-\kappa_u)} N_\tau +  \frac{1}{1-\kappa_y} \Theta_\tau. \numberthis \label{eq:ys-theta}
\]
We now fix $y^s_{-1} = y^s_0$. Then, \eqref{eq:ys-theta} immediately implies
\begin{align*}
	&\sum_{t=0}^\tau \norm{y^s_{t} -y^s_{t-1}} \leq 2\sum_{t=0}^\tau \norm{y^s_{t} -\theta_{t}} + \sum_{t=0}^\tau \norm{\theta_{t} -\theta_{t-1}} \\
	{\refleq{\eqref{eq:ys-theta}}} &\frac{c_1\kappa_y\sqrt{\tilde\mu}}{(1{-}\kappa_y)(1{-}\kappa_u)} N_\tau {+}  \frac{3{-}\kappa_y}{1{-}\kappa_y} \Theta_\tau. \numberthis \label{eq:ys-ys}
\end{align*}
Moreover, we define $S_0 = I_m - S_u^\dagger S_u$ and have
\begin{align*}
	&\sum_{t=0}^\tau \norm{u_t^s - v_t} \leq \sum_{t=0}^\tau \norm{u_t^s - \eta_t} + \sum_{t=0}^\tau \norm{v_t - \eta_t} \\
	\refleq{\eqref{eq:SSInput},\eqref{eq:Lemma1}} &\!\left( \norm{S_0} + 1\right) \! \sum_{t=0}^\tau \norm{v_t - \eta_t} + \norm{S_u^\dagger S_y} \sum_{t=0}^\tau \norm{y^s_t - \theta_t} \\
	\refleq{\eqref{eq:vt-etat},\eqref{eq:ys-theta}} & \! \left( \! \frac{\norm{S_0} {+} 1}{1{-}\kappa_u} {+} \frac{c_1\kappa_y\sqrt{\tilde\mu}\norm{S_u^\dagger S_y}}{2(1{-}\kappa_y)(1{-}\kappa_u)}\!\right) \! N_\tau {+} \frac{\norm{S_u^\dagger S_y}}{1{-}\kappa_y} \Theta_\tau. \numberthis \label{eq:us-v}
\end{align*}
Let $\bar H_\alpha = Y^{\tilde\mu:\tilde\mu+n-1}H_\alpha^\dagger$. Then we get
\begin{align*}
	&\sum_{t=0}^\tau \norm{g_t} = \sum_{t=0}^\tau \norm{\begin{bmatrix} 1_{n+1} \otimes u_t^s - U^{\tilde\mu:\tilde\mu+n}(\alpha_t+\omega_t) \\ 1_{n} \otimes y^s_t - Y^{\tilde\mu:\tilde\mu+n-1} (\alpha_t+\omega_t) \end{bmatrix} } \\
	&\refeq{\eqref{eq:defalpha},\eqref{eq:RecursivePredictions}} \sum_{t=0}^\tau \norm{\begin{bmatrix} 1_{n+1} \otimes (u_t^s - v_t - u^s_{t-1} + v_{t-1}) \\ 1_{n} \otimes (y^s_t {-} y^s_{t-1}) + Y^{\tilde\mu:\tilde\mu +  n - 1} ( \omega_t{-}\omega_{t-1}) \end{bmatrix}} \\
	\begin{split} 
		 &\refleq{\eqref{eq:defomega}} \sqrt{n+1} \sum_{t=0}^\tau \Big( \norm{u^s_t - u^s_{t-1}} + \norm{v_{t+1} - v_{t}} \Big) \\
		 &+ \sum_{t=0}^\tau \left( \sqrt{n} \norm{y^s_t - y^s_{t-1}} + \sqrt{\tilde\mu} \norm{\bar H_\alpha}  \norm{v_{t+1} - v_{t}} \right)
	\end{split} \\
	\begin{split}
		&\refleq{\eqref{eq:SSInput}} \sqrt{\tilde\mu} \bigg(  \norm{S_0}  {+} \norm{\bar H_\alpha} {+} 1 \bigg) \sum_{t=0}^\tau \norm{v_{t+1} {-} v_{t}} \\
		&+ \sqrt{\tilde\mu} \bigg( \norm{S_u^\dagger S_y} + 1 \bigg) \sum_{t=0}^\tau \norm{y^s_t - y^s_{t-1}}
	\end{split}\\
	&\refleq{\eqref{eq:vt-vt},\eqref{eq:ys-ys}} \hspace{-3pt} \frac{\sqrt{\tilde\mu}\tilde c_2(1{-}\kappa_y){+}\kappa_y \tilde\mu c_1\tilde c_3}{(1-\kappa_u)(1-\kappa_y)} \! N_\tau {+} \sqrt{\tilde\mu}\frac{3-\kappa_y}{1-\kappa_y}\tilde c_3 \Theta_\tau,
\end{align*}
where $\tilde c_2 := l_u (\norm{S_0} + \norm{\bar H_\alpha} + 1)$ and $\tilde c_3 := (\norm{S_u^\dagger S_y} + 1)$.
Let $\bar Q= \left( I - \left( Q \left( I - H_\beta^\dagger H_\beta \right) \right)^\dagger Q \right) H_\beta^\dagger $, $c_2 = \norm{\bar Q}\tilde c_2$, and $c_3 = \norm{\bar Q}\tilde c_3$. Since $\beta_t = \bar Q g_t$, this immediately implies
\begin{equation}
	 \sum_{t=0}^\tau \norm{\beta_t} \leq \hspace{-3pt} \frac{\sqrt{\tilde\mu} c_2(1{-}\kappa_y){+}\kappa_y \tilde\mu c_1c_3}{(1-\kappa_u)(1-\kappa_y)} \! N_\tau {+} \sqrt{\tilde\mu}\frac{3-\kappa_y}{1-\kappa_y}c_3 \Theta_\tau. \label{eq:boundbeta}
\end{equation}

Next, repeatedly inserting \eqref{eq:PredInputs} yields
\begin{align*}
	\hat u_t(1) 
	&= u^s_{t-\mu-1} - v_{t-\mu-1} + \sum_{i=0}^{\mu-1} U^{n+1+i} \beta_{t-i}.
\end{align*}
Finally, we can combine the previous results together with $v_t = u^s_t$ and $\beta_t = 0$ if $t<0$ to get	
\begin{align*}
	&\sum_{t=0}^{\tau} \norm{u_t-\eta_t} = \sum_{t=0}^{\tau} \norm{\hat u_t(1)+v_t-\eta_t} \\
	\refleq{\eqref{eq:vt-etat}} & \frac{1}{1{-}\kappa_u} N_\tau {+} \sum_{t=0}^{\tau} \left( \norm{\sum_{i=0}^{\mu-1} \hspace{-3pt} U^{n+1+i} \beta_{t-i}} {+} \norm{u^s_{t}{-}v_{t}} \right) \\
	&\refleq{\eqref{eq:us-v},\eqref{eq:boundbeta}} \hspace{-2pt} \frac{1}{1{-}\kappa_u} \! \left( \! \norm{S_0} {+} 2 {+} \frac{c_1\kappa_y\sqrt{\tilde\mu}\norm{S_u^\dagger S_y}}{2(1{-}\kappa_y)}\!\right) \! N_\tau \\
	&+ \mu \norm{U^{n+1:\tilde\mu}} \frac{\sqrt{\tilde\mu} c_2(1-\kappa_y)+\kappa_y \tilde\mu c_1c_3}{(1-\kappa_u)(1-\kappa_y)}  N_\tau \numberthis \label{eq:u-eta} \\
	&+ \left( \frac{\norm{S_u^\dagger S_y}}{1-\kappa_y} + \mu \norm{U^{n+1:\tilde\mu}} \sqrt{\tilde\mu}\frac{3-\kappa_y}{1-\kappa_y}c_3 \right) \Theta_\tau. 
\end{align*}

\subsubsection*{Part 2 (Predicted Regret)}
Since $\theta_{-1} = \hat y_0^\mu$, we get
\begin{align*}
	&\sum_{t=0}^\tau \norm{\hat y^\mu_{t} - \theta_{t-1}} = \sum_{t=0}^\tau \norm{\hat y^\mu_{t+1} - \theta_{t}} \\
	\refleq{\eqref{eq:OutputPrediction},\eqref{eq:RecursivePredictions}} &\sum_{t=0}^\tau \norm{y^s_t - \theta_{t-1}} + \sum_{t=0}^\tau \norm{Y^{\tilde\mu} (\omega_{t+1} - \omega_t)} + \Theta_\tau \\
	\begin{split}
		\refleq{\eqref{eq:ContractionGD},\eqref{eq:vt-vt}} &\kappa_y \sum_{t=0}^\tau \norm{\hat y_{t}^\mu - \theta_{t-1}} + \norm{Y^{\tilde\mu} H_\alpha^\dagger} \frac{\sqrt{\tilde\mu} l_u}{1-\kappa_u} N_\tau + \Theta_\tau
	\end{split}
\end{align*}
Due to $1-\kappa_y > 0$, we can rearrange to get the result
\[
	\sum_{t=0}^\tau \norm{\hat y^\mu_{t} - \theta_{t-1}} \leq \frac{1}{1{-}\kappa_y} \Theta_\tau + \frac{\sqrt{\tilde\mu}c_1/2}{(1{-}\kappa_u)(1{-}\kappa_y)} N_\tau. \numberthis \label{eq:yhat-theta}
\]

\subsubsection*{Part 3 (Prediction Error)}
By Willems' Lemma, it is possible to write the real output $y_{t+\mu}$ at time $t+\mu$ as $y_{t+\mu} = Y^{\tilde\mu} \bar \alpha_t$, where $\bar \alpha_t$ satisfies
\[
\bar \alpha_t = H_{\alpha}^\dagger \begin{bmatrix} u_{[t-n:t+\mu]} \\ 1_{n} \otimes (u^s_{t-1} - v_{t-1} + v_t) \\ y_{[t-n:t-1]} \end{bmatrix}.
\]
Moreover, the real inputs $u$ can be expressed as
\begin{align*}
	u_{t+i} &= (\sigma \hat u_{t-1})_{i{+}1} + \sum_{j=0}^{i} \Big( U^{n+i+1-j} \beta_{t+j} \Big) + v_{t+i}, \\
	u_{t+\mu} &= u^s_{t-1}-v_{t-1} + \sum_{j=0}^{\mu} \Big( U^{\tilde\mu-j} \beta_{t+j} \Big) + v_{t+\mu},
\end{align*}
where $0 \leq i \leq \mu-1$. Then, we obtain the bound as follows
\begin{align*}
	&\sum_{t=0}^\tau \norm{y_{t+\mu} - \hat y_{t+\mu}} = \sum_{t=0}^\tau \norm{Y^{\tilde\mu}(\bar \alpha_t - \alpha_t - \omega_t)} \\
	&\leq \norm{Y^{\tilde\mu}H_\alpha^\dagger} \sum_{t=0}^\tau \norm{
		\begin{bmatrix}
			U^{n+1} \beta_t \\
			\dots \\
			v_{t+\mu} - v_t + \sum_{j=0}^{\mu} U^{\tilde\mu-j} \beta_{t+j}
	\end{bmatrix}} \\
	&\leq \norm{Y^{\tilde\mu}H_\alpha^\dagger} \sum_{t=0}^{\tau}  \sum_{i=0}^\mu  \bigg( \norm{U^{n+1:\tilde\mu}\beta_{t+i}} + \norm{v_{t+i}{-}v_t}  \bigg)\\
	&\refleq{\eqref{eq:vti-vt}} c_4 (\mu{+}1) \sum_{t=0}^{\tau+\mu} \norm{\beta_t} {+} \frac{\mu^2 c_1}{2(1{-}\kappa_u)} N_{\tau+\mu} \\
	&\refleq{\eqref{eq:boundbeta}} \! \left(c_4(\mu{+}1) \frac{\sqrt{\tilde\mu} c_2(1{-}\kappa_y){+}\kappa_y \tilde\mu c_1c_3}{(1-\kappa_u)(1-\kappa_y)} + \frac{\mu^2 c_1}{2(1{-}\kappa_u)} \right) N_{\tau+\mu} \\
	&+ c_4 (\mu+1) \sqrt{\tilde\mu}\frac{3-\kappa_y}{1-\kappa_y} c_3 \Theta_{\tau+\mu},\numberthis \label{eq:y-yhat}
\end{align*}
with $c_4 = \norm{Y^{\tilde\mu}H_\alpha^\dagger} \norm{U^{n+1:\tilde\mu}}$. The theorem follows from inserting \eqref{eq:thetamu-theta}, \eqref{eq:u-eta}, \eqref{eq:yhat-theta}, and \eqref{eq:y-yhat} into \eqref{eq:RegretUnfinished}. \qed


\bibliographystyle{IEEEtran}
\bibliography{IEEEabrv,IEEEexample}


\end{document}